\documentclass[12pt,leqno]{article}
\usepackage{amsmath,amssymb,amsfonts}
\usepackage{eucal}

\usepackage {epsfig}
\newcommand\comment[1]{}                %  Silent version.
\renewcommand\comment[1]{\emph{[#1]}}           %  Comment revealed.

\numberwithin{equation}{section}

\newtheorem{theorem}{Theorem}

\newtheorem{definition}[theorem]{Definition}

\newtheorem{lemma}[theorem]{Lemma}

\newtheorem{problem}[theorem]{Problem}
\newtheorem{proposition}[theorem]{Proposition}
\newtheorem{remark}[theorem]{Remark}

\newenvironment{proof}[1][Proof]{\textbf{#1.} }{\ \rule{0.5em}{0.5em}}

\renewcommand{\phi}{\varphi}                 % Personal preferences.
\renewcommand{\epsilon}{\varepsilon}

\newcommand\const{\operatorname{const}}

\newcommand\Star{\operatorname{Star}}
\newcommand\conv{\operatorname{conv}}

\newcommand\PStar{\operatorname{PStar}}

\newcommand\GLB{\operatorname{GLB}}

\newcommand\x{\mathbf{x}}

\newcommand\q{\mathbf{q}}

\newcommand\p{\mathbf{p}}
\newcommand\oo{\mathbf{o}}

\newcommand\R{\mathbb{R}}
\newcommand\X{\mathbb{X}}
\newcommand\HH{\mathbb{H}}
\newcommand\SSS{\mathbb{S}}

\newcommand\cM{\mathcal{M}}

\begin{document}

\title{On locally convex PL-manifolds \\and fast verification of convexity}
\author{Konstantin Rybnikov \\
\date{\today}
krybniko@cs.uml.edu\\
http://faculty.uml.edu/krybnikov}
\maketitle

\centerline{Short Version}
\medskip

\begin{abstract}
We show that a PL-realization of a closed connected manifold of
dimension $n-1$ in $\mathbb{R}^n\:(n \ge 3)$ is the boundary of a
convex polyhedron if and only if the interior of  each
$(n-3)$-face has a point, which has a neighborhood lying on the
boundary of a convex $n$-dimensional body. This result is derived from a
generalization of Van Heijenoort's theorem on locally convex
manifolds to the spherical case. Our convexity criterion
for PL-manifolds implies an easy polynomial-time algorithm for
checking convexity of a given  PL-surface in $\mathbb{R}^n$.
\end{abstract}

There is a number of theorems that infer global convexity from local convexity. The oldest one belongs to Jacque Hadamard (1897) and asserts that any compact smooth surface embedded in $\R^3$, with strictly positive Gaussian curvature, is the boundary of a convex body.  Local convexity can be defined in many different ways (see van Heijenoort (1952) for a survey).  We will use Bouligand's (1932) notion of local convexity.  In this definition a surface $M$ in the affine space $\mathbb{R}^n$ is called locally convex at point ${\bf p}$ if ${\bf p}$ has a neighborhood which lies on the boundary of a convex $n$-dimensional body $K_{\bf p}$; if  $K_{\bf p}\backslash{\bf p}$ lies in an open half-space defined by a hyperplane containing  ${\bf p}$, $M$ is called strictly convex at ${\bf p}$.

This paper is mainly devoted to local convexity of piecewise-linear (PL) surfaces, in particular, polytopes. A PL-surface in $\mathbb{R}^n$ is a pair $M=({\cal M},r)$, where ${\cal M}$ is a topological manifold with a fixed cell-partition and $r$ is a continuous {\it realization} map from ${\cal M}$ to  $\mathbb{R}^n$ that satisfies the following conditions:
\par \noindent 1) $r$ is a bijection on the closures of all cells of  ${\cal M}$
\par \noindent 1) for each $k$-cell $C$ of $\cM$ the image $r(C)$ lies on a $k$-dimensional affine subspace of $\mathbb{R}^n$; $r(C)$ is then called a $k$-face of $M$.
\par Thus, $r$ need not be an immersion, but its restriction to the closure of any cell of  ${\cal M}$ must be. By a fixed cell-partition of ${\cal M}$ we mean that ${\cal M}$ has a structure of a CW-complex where all gluing mappings are homeomorphisms (such complexes are called \emph{regular} by J.H.C. Whitehead).  \emph{All cells and faces are assumed to be open.} We will also call $M=({\cal M},r)$ a PL-realization of ${\cal M}$ in $\mathbb{R}^n$.

\begin{definition}
\emph{We say that $M=({\cal M},r)$ is the boundary of a convex body $P$ if $r$ is a homeomorphism between  ${\cal M}$ and ${\partial P}$.} 
\end{definition}

Hence, we exclude the cases when $r(\cM)$ coincides with the boundary of a convex set, but $r$ is not injective. Of course, the algorithmic and topological  sides of this case are rather important for  computational geometry and we will consider them in further works. Notice that for $n>2$ a closed $(n-1)$-manifold $\cM$ cannot be immersed into $\R^n$ by a 
non-injective map $r$ so that $r(\cM)$ is the boundary of a convex set, since any covering space of a simply connected manifold must be simply connected. However, such immersions are possible in the hyperbolic space $\HH^n$.

Our main theorem asserts that any closed PL-surface $M$ immersed in $\mathbb{R}^n \:(n \ge 3)$ with at least one point of strict convexity, and such that each $(n-3)$-cell has a point at which $M$ is locally convex, is convex. Notice that if the last condition holds for some point on an $(n-3)$-face, it holds for all points of this face.  

This theorem implies a test for global convexity of PL-surfaces: check local convexity on each of the $(n-3)$-faces. Notice that if for all $k$ and  every $k$-face there is an $(n-k-1)$-sphere, lying in a complementary subspace and centered at some point of $F$, such that $\SSS \cap F$ is a convex surface, then  $r$ is an immersion. The algorithm implicitly checks if a given realization is an immersion, and reports "not convex" if it is not. 

 The pseudo-code for the algorithm is given in this article. The complexity of this test depends on the way the surface is given as input data. Assuming we are given the coordinates of the vertices and the poset of faces of dimensions $n-1$, $n-2$, $n-3$ and $0$, OR, the equations of the facets, and the poset of faces dimensions $n-1$, $n-2$, and $n-3$, the complexity of the algorithm for a general closed PL-manifold is $O(f_{n-3,n-2})=O(f_{n-3,n-1})$, where $f_{k,l}$ is the number of incidences between cells of dimension $k$ and $l$.
 If the vertices of the manifold are assumed to be in a sufficiently general position, then the dimension of the space does not affect the complexity at all. Another advantage of this algorithm is that it consist of $f_{n-3}$ independent subroutines corresponding to the $(n-3)$-faces, each with complexity not exceeding $O$ in the number of $(n-1)$-cells incident to the $(n-3)$-face.

 The complexity of our algorithm is asymptotically equal to the complexity of algorithms suggested by Devillers et al (1998) and Mehlhorn et al (1999) for simplicial 2-dimensional surfaces; for $n>3$ our algorithm is asymptotically faster than theirs. These authors verify convexity not by checking it locally at $(n-3)$-faces, but by different, rather global methods (their notion of local convexity is, in fact, a global notion). Devillers et al (1998) and Mehlhorn et al (1999) make much stronger initial assumptions about the input, such as the orientability of the input surface;  they  also presume that for each $(n-1)$-face of the surface an external normal is given, and that the directions of these normals define an orientation of the surface. Then they call the surface locally convex at an $(n-2)$-face $F$ if the angle between the normals of two $(n-1)$-faces adjacent to $F$ is obtuse. Of course this notion of ``local'' convexity is not local.

The main theorem is deduced from a direct generalization of van Heijenoort's theorem to the spherical case.
Van Heijenoort's theorem asserts that an immersion in $\mathbb{R}^n$ of any closed connected manifold ${\cal M}$, which is locally convex at all points, strictly locally convex in at least one point, and is complete with respect to the metric induced by $r$, is the boundary of a convex $d$-dimensional set. Van Heijenoort (1952) noticed that for $n=3$ his theorem immidiately follows from four theorems contained in Alexandrov (1948); however, acoording to van Heijenoort, Alexandrov's methods do not extend to $n>3$ and his approach is technically more complicated. We show that this theorem also holds for  spheres, but not for the hyperbolic space. While all notions of affine convexity  can be obviously generalized to the hyperbolic space, there are two possible generalizations in the spherical case; neither of these generalizations is perfect. The main question is whether we want all geodesics joining points of a convex set $S$ to be contained in $S$, or at least one. In the first case subspheres are not convex, in the second case two convex sets can intersect by a non-convex set. The latter problem can be solved by requiring a convex set to be open, but this is not very convenient, since again, it excludes subspheres. We call a set in $\X^n$ convex if for any two points $\p,\q \in S$, there is some geodesic $[\p,\q] \subset S$. 
\begin{proposition}
If the intersection $I$ of two convex sets in $\SSS^n$, $n>0$ is not convex, then $I$ contains two opposite points.
\end{proposition}

To have unified terminology we will call subspheres subspaces.

Besides the algorithmic implications, our generalization implies that any $(n-3)$-simple PL-surface in $\mathbb{R}^n$  with convex facets is the boundary of a convex polyhedron.

\section{van Heijenoort-Alexandrov's Theorem \\for Spaces of Constant Curvature}
Throughout the paper $\X^n$ denotes $\R^n$, $\SSS^n$, or $\HH^n$.
Following the original proof of van Hejeenoort's, we will now show that his theorem holds in a somewhat stronger form for $\SSS^n$ for $n>2$. van Hejeenoort's theorem does not hold for unbounded surfaces in $\HH^n$. We will give three different kinds of counterexamples and pose a conjecture about simply connected locally compact embeddings of manifolds in $\HH^n$.

Imagine a "convex strip" in 3D which is bent in the form of handwritten  $\varphi$ so that it self-intersects itself, but not locally. Consider the intersection of this strip with a ball of appropriate radius so that the self-intersection of the strip happens to be inside the ball, and the boundary of the strip  outside. Regarding the interior of the ball as Klein's model of $\HH^3$ we conclude that the constructed surface is strictly locally convex at all points and has a complete metric induced by the immersion into the hyperbolic space. This gives an example of an\emph{ immersion of a  simply connected manifold} into $\HH^n$ which does not bound a convex surface. Notice that in this counterexample the surface self-intersects itself. 

Consider the (affine) product of a non-convex quadrilateral, lying inside a unit sphere centered at the origin, and a line in $\R^n$. The result is a non-convex polyhedral cylindrical surface $P$. Pick a point $\p$ inside the sphere, but outside the cylinder, whose vertical projection on the cylinder is the affine center of one of its facets. Replace this facet of $P$ with the cone over $F \cap \{\x | \|\x\| \le 1\}$ The part of the resulting polyhedral surface, that lies inside the sphere of unit radius, is indeed a PL-surface \emph{embedded} in $\HH^n$ (in Klein's model). The surface if locally convex at every point and strictly convex at $\p$. However, it is not the boundary of a convex body. Notice that this surface is \emph{not simply connected}. 

Consider a locally convex spiral, embedded in the $yz$-plane with two limiting sets: circle $\{\p | x=0, y^2+z^2=1  \}$ and the origin. That is this spiral coils around the origin and also around (from inside) the circle. Let $M$ be the double cone over this spiral with apexes at $(1,0,0)$ and $(-1,0,0$, intersected with the unit ball $\{\x | \|\x\| < 1\}$. This \emph{non-convex } surface is obviously \emph{simply connected, embedded,} locally convex at every point, and strictly convex at all point of the spiral. 

A locally compact realization $r$ of $\cM$ is a realization such that for any compact subset $C$ of $\X^n$ $C \cap r(\cM)$ is compact. The question remains:
\begin{problem} Is it true that any locally compact embedding of a simply connected surface in $\HH^n$ is convex?  
\end{problem}
It remains an open question whether van Heijenoort's (1952) criterion works for \emph{embedded} \emph{unbounded} surfaces in $\HH^n$. We conjecture that it is, indeed, the case. The proof of the main theorem makes use of quite a number of technical propositions and lemmas. The proofs of these statements for $\R^n$ by most part can be directly repeated for $\X^n$, but in some situations  extra care is needed. If the reader is referred to van Heijenoort's paper for the proof, it means that the original proof works without any changes.

\underline{Notation:} The calligraphic font is used for sets in the abstract topological manifold. The regular mathematics font is used for the images of these sets in $\X^n$, a space of constant curvature. The interior of a set $S$ is denoted by $(S)$, while the closure by $[S]$. The boundary of S is denoted by $\partial S$. Since this paper is best read together with van Heijenoort's (1952) paper, we would like to explain the differences between his and our notations. van Heijenoort denotes a subset in the abstract manifold $\overline{M}$ by $\overline{S}$, while denoting its image in $\mathbb{R}^n$ by $S$;  an interior of a set $S$ in $\mathbb{R}^n$ is denoted in his paper by $\dot{S}$.

The immersion $r$ induces a metric on ${\cal M}$ by \[ d(p,q)=\GLB\limits_{arc(p,q) \subset {\cal M}}\{|r(arc(p,q))|\} \]
where \emph{$\GLB$} stands for the greatest lower bound, and $|r(arc(p,q))|$ — for the length of an arc joining $\p$ and $\q$ on $M$, which is the $r$-image of an arc joining these points on ${\cal M}$. We will call this metric $r$-metric.

\begin{lemma}\label{arcwise} (van Heijenoort) Any two points of ${\cal M}$ can be connected by an arc of a finite length. Thus  ${\cal M}$  is not only connected, but also arcwise connected.
\end{lemma}

\begin{lemma} (van Heijenoort) The metric topology defined by the $r$-metric is equivalent to the original topology on ${\cal M}$.
\end{lemma}

\begin{lemma} (van Heijenoort) $r(S)$ is closed in $\X^n$ for any closed subset $S$ of ${\cal M}$.
\end{lemma}

\begin{lemma} (van Heijenoort) If on a bounded (in $r$-metric)  closed  subset $S \subset {\cal M}$ mapping $r$ is one-to-one, then $r$ is a homeomorphism between  $S$ and $r(S)$.
\end{lemma}

The proofs of the last two lemmas have been omitted in van Heijenoort (1952), but they are well known in topology.

\begin{theorem}\label{main}
 Let $\X^n$ ($n>2$) be a Euclidean, spherical, or hyperbolic space. Let $M=({\cal M},r)$ be an immersion of an $(n-1)$-manifold ${\cal M}$ in $\X^n$, such that $r(\cM)$ is bounded in $\X^n$. Suppose that $M=({\cal M},r)$ satisfies the following conditions:

\noindent  1) ${\cal M}$ is complete with respect to the metric induced on ${\cal M}$ by the immersion $r$,

\noindent  2) ${\cal M}$ is connected,

\noindent  2) $M$ is locally convex at each point,

\noindent  4) $M$ is strictly convex in at least one point,

Then $r$ is a homeomorphism from ${\cal M}$ onto the boundary of a compact convex body. 
\end{theorem}
\begin{proof}
Notice that our theorem for $\X^n=\HH^n$ directly follows from van Heijenoort's proof of the Euclidean case. Any immersion of $\cM$ into $\HH^n$ can be regarded as an immersion into  the interior of a unit ball with a hyperbolic metric, according to Klein's model. If conditions 1)-4) are satisfied for the hyperbolic metric, they are satisfied for the Euclidean metric on this ball. Geodesics in Klein's model are straight line segments and, therefore, for a bounded closed surfaces in $\HH^n$, that satisfies the conditions of the theorem, the convexity follows from the Euclidean version of this theorem.

The original Van Heijenoort's proof is based on the notion of convex part. A {\it convex part} of  $M$, centered at a point of strict convexity $\oo=r(o)$, $o \in \cM$, is an open connected subset $C$ of  $r({\cal M})$ that contains   $\oo$  and such that: (1) $\partial C= H \cap r({\cal M})$, where $H$ is a hyperplane in $\X^n$, not passing through $\mathbf{o}$, (2) $C$ lies on the boundary of a closed convex body $K_C$ bounded by $C$ and $H$. We call $H \cap K_C$ the lid of the convex part. Let $H_0$ be a supporting hyperplane at $\oo$. We call the \emph{open} half-space defined by $H_0$, where the convex part lies, the \emph{positive half-space} and denote it by $H^+_0$. We call the $r$-preimage of a convex part $C$ in ${\cal M}$ an abstract convex part,  and denote it by ${\cal C}$. In van Heijenoort's paper $H$ is required to be parallel to the supporting hyperplane $H_0$ of $r({\cal M})$  at $\oo$, but this is not essential. In fact, we just need a family of hyperplanes such that: (1) they do not intersect in the positive half-space, (2) the intersections of these hyperplanes with the positive half-space form a partition of the positive half-space, (3) all these hyperplanes are orthogonal to a line $l$, passing through $\oo$. Let us call such a family a fiber bundle $\{H_z\}_{(l,H_0)}$ of the positive half-space defined by $l$ and $H_0$. (In the case of $\X^n=\R^n$ it is a vector bundle.) In fact, it is not necessary to assume that $l$ passes through $\oo$, but this assumption simplifies our proofs. Here $z>0$ denotes the distance, \emph{along the line $l$}, between the hyperplane $H_z$ in this family  and $H_0$. We will call $z$ the height of $H_z$.
\begin{proposition}
A convex part exists.
\end{proposition}
\begin{proof}van Heijenoort's proof works for $\X^n$, $n>2$, without changes.
\end{proof}

Denote by $\zeta$ the least upper bound of the set of heights of the lids of convex parts centered at $\oo$ and defined by some fixed fiber bundle $\{H_z\}_{(l,H_o)}$. Since $r(\cM)$ is bounded, then  $\zeta < \infty$.

Consider the union $G$ of all convex parts, centered at $\oo$. We want to prove that this union is also a convex part. Let us depart for a short while (this paragraph) from the assumption that $r(\cM)$ is bounded.  $G$ may only be unbounded in the hyperbolic and Euclidean cases.  As shown by van Heijennort (1952), if $\X^n=\R^n$ and $\zeta < \infty$,  $G$ must be bounded  even when  $r(\cM)$ is allowed to be unbounded. If $\X^n=\HH^n$ and $\zeta < \infty$,  $G$ can be unbounded, and this is precisely the reason why van Heijenoort's theorem does not hold for unbounded surfaces in hyperbolic spaces.

%The following is a lie:
%If  $G$ is unbounded, then, following, word by word, van %Heijenoort's argument for the Euclidean case (part I of %Section 3, 1952), we conclude that $r({\cal M})$ is %homeomorphic to a hyperplane in $\X^n$. His proof is equally %good for the hyperbolic  case.

Since in this theorem $r(\cM)$ is assumed to be bounded, $G$ is bounded.
Let us presume from now on that $\X^n=\SSS^n$ (the case of $\X^n=\HH^n$ is considered in the beginning of the proof).
 $\partial G$  belongs to the hyperplane $H_{\zeta}$ and is equal to $H_{\zeta} \cap M$. $\partial G$ bounds a closed bounded convex set $D$ in $H_{\zeta}$. Two mutually excluding cases are possible.

Case 1: $\dim D<n-1$. Then, following the argument of van Heljenoort (Part 2: pages 239-230, Part 5: page 241, Part 3: II on page 231), we conclude that $G \cup D$ is the homeomorphic pre-image of an $(n-1)$-sphere ${\cal G} \cup {\cal D} \subset {\cal M}$. Since ${\cal M}$ is connected, ${\cal G} \cup {\cal D} = {\cal M}$, and ${\cal M}$ is a convex surface.

Case 2: $\dim D=n-1$. The following lemma is a key part of the proof of the main theorem. Roughly speaking, it asserts that 
if the lid of a convex part is of co-dimension 1, then either this convex part is a subset of a bigger convex part, or this convex part, together with the lid, is homeomorphic to $\cM$ via mapping $r$.

\begin{lemma}\label{alternative} Suppose $\X^n=\SSS^n$. Let $C$ be a convex part centered at a point $\oo$ and defined by a hyperplane $H_z$ from a fiber bundle $\{H_z\}_{(l,H_0)}$. Suppose  $B=\partial C$ is the boundary of an $(n-1)$-dimensional closed convex set $S$ in $H_z$.  Either $S$ is the $r$-image of an $(n-1)$-disk ${\cal S}$ in ${\cal M}$ and ${\cal M}={\cal C} \cup {\cal S}$, where ${\cal C}=r(C)$, or $C$ is a proper subset of a larger convex part, defined by the same fiber bundle $\{H_z\}_{(l,H_0)}$.
\end{lemma}

\begin{proof} Using a perturbation argument, we will prove this lemma by reducing the spherical case to the Euclidean one.
Since $S$ is $(n-1)$-dimensional and belongs to one of the hyperplanes in the fiber bundle $\{H_z\}_{(i,H_0)}$, $[\conv C] \cap H_0$ is either empty or $(n-1)$-dimensional. If it is non-empty, $[\conv C] \cap H_0$ must have a point other than $\oo$ and its opposite. The closure of a convex set in $\X^n$ is convex.  Since $[\conv C ]$ is convex, if it contains a point $\p$ of $H_0$ other than $\oo$ and its opposite, it contains some  geodesic segment $[\oo \p]$ \emph{lying in} $H_0$. Since $\oo$ is a point of strict convexity, there is a neighborhood of $\oo$ on $\oo \p$ all whose points, except for $\oo$,  are not  points of $[\conv C]$, which contradicts to the choice of  $[\oo \p]$. 

  So, $[\conv C] \cap H_0$ is definitely empty.  Since, by Lemma \ref{arcwise}, $r(M)$ is arcwise connected,
all of $[\conv C]$, except for the point $\oo$, lies in the positive subspace.  Therefore, there is a hyperplane $H$ in $\SSS^n$ such that $[C]$ lies in an open halfspace $H_+$ defined by $H$. We can regard $\SSS^n$ as a standard sphere in $\R^{n+1}$. $H$ defines a hyperplane in $\R^{n+1}$. Consider an $n$-dimensional plane $E_n$ in $\R^{n+1}$ parallel to this hyperplane and not passing through the origin. Central projection $r_1$ of $M \cap H_+$ on $E_n$ obviously induces an immersion $r_1r$ of a submanifold $\cM^{\prime}$ of $\cM$ into $E_n$.

This submanifold $\cM^{\prime}$ is defined as the maximal  arcwise connected open subset of $\cM$ such that (1) all points of this subset are mapped by $r$ to $H_+$, and  (2) it contains $o$. It is obviously a manifold. Let us prove that it exists. Consider the union of all open arcwise connected subsets that contain $o$. It is open and is acrwise connected, since it contains $o$. Let $M^{\prime}=(\cM^{\prime}, r_1r)$.

The immersion $r_1r$ obviously satisfies Conditions 2-4 of the main theorem \ref{main}. $r_1r$ defines a metric on $\cM^{\prime}$. Any Cauchy sequence on $\cM^{\prime}$ under this metric  is also a Cauchy sequence on $\cM$ under the metric induced by $r$. Therefore $\cM^{\prime}$ is complete and satisfies the conditions of the main Theorem \ref{main}. The central projection on $E_n$ maps a spherical convex part of $M$ on a Euclidean convex part of $M^{\prime}$; it also maps the fiber bundle $\{H_z\}_{(l,H_0)}$ to a fiber bundle in the Euclidean $n$-plane $E_n$. van Heijenoort (1952) proved Lemma \ref{alternative}  for the Euclidean case. Therefore, either ${\cal M}={\cal C} \cup {\cal S}$, or $C$ is a proper subset of a larger convex part centered at $\oo$, and defined by the same fiber bundle $\{H_z\}_{(l,H_0)}$. 
\end{proof}

The second alternative ($C$ is a subset of a larger convex part) is obviously excluded, since $C$ is the convex part corresponding to the height which is the least upper bound of all possible heights of convex parts. Therefore in Case 2 ${\cal M}$ is the boundary of a convex body which consists of a maximal convex part and a convex $(n-1)$-disk, lying in the hyperplane $H_{\zeta}$.

\end{proof}

\section{Locally convex PL-surfaces}
\begin{theorem}\label{strict} Let $r$ be a realization map from a \emph{compact} connected manifold ${\cal M}$ of dimension $n-1$ into $\mathbb{X}^n$ ($n>2$) such that $\cM$ is complete with respect to the $r$-metric.  Suppose that $M=(\cM,r)$ is locally convex at all points. Then  $M=(\cM,r)$ is either strictly locally convex in at least one point, or is a spherical hyper-surface of the form $\mathbb{S}^n \cap \partial C$, where $C$ is a convex cone in $\mathbb{R}^{n+1}$, whose face of the smallest dimension contains the origin (in particular, $C$ may be a hyperplane in $\mathbb{R}^{n+1}$). 
\end{theorem}
\begin{proof} The proof of this rather long and technical theorem will be included in the full length paper of Rybnikov (200X).
\end{proof}

\begin{theorem}\label{PL-case} Let $r$ be a realization map from a closed connected $n$-dimensional manifold ${\cal M}$, with a regular CW-decomposition,  in $\mathbb{R}^n$ or $\SSS^n$ ($n>2$) such that on the closure of each cell $C$ of ${\cal M}$ map $r$ is one-to-one and $r(C)$ lies on a subspace of dimension equal to $\dim C$.  Suppose that $r({\cal M})$ is strictly locally convex in at least one point. The surface $r({\cal M})$ is the boundary of a convex polyhedron if and only if each $(n-3)$-face has a point with an  $M$-neighborhood which lies on the boundary of a convex $n$-dimensional set.
\end{theorem}
\begin{proof}
$M$ is locally convex at all points of its $(n-3)$-cells. Suppose we have shown that $M$ is locally convex at each $k$-face, $0<k \le n-3$. Consider a $(k-1)$-face $F$.  Consider the intersection of $\Star(F)$ with a sufficiently small $(n-k)$-sphere ${\mathbb S}$ centered at some  point $\p$ of $F$ and lying in a subspace complimentary to $F$. $M$ is locally convex at $F$ if and only if the hypersurface  ${\mathbb S} \cap \Star(F)$ on the sphere  ${\mathbb S}$ is convex. Since   $M$ is locally convex at each $k$-face, ${\mathbb S} \cap \Star(F)$ this hypersurface is locally convex at each vertex. By Theorem \ref{strict} ${\mathbb S} \cap \Star(F)$ is either the intersection of the boundary of a convex cone with ${\mathbb S}$ or has a point of strict convexity on the sphere ${\mathbb S}$.  In the latter case the spherical generalization of van Heijenoort's theorem implies that ${\mathbb S} \cap \Star(F)$ is convex. Thus $M$ is convex at $\p$ and therefore at all points of $F$.

This induction argument shows that $M$ must be locally convex at all vertices. If is locally convex at all vertices, it is locally convex at all points. We assumed that $M$ had a point of strict convexity. The metric induced by $r$ is indeed complete. By van Heijenoort's theorem and Theorem \ref{main}, $M$ is the boundary of a convex polyhedron.
\end{proof}

\section{New Algorithm for Checking Global \\Convexity of PL-surfaces}
\underline{Idea:} check convexity for the star of each $(n-3)$-cell of $M$.

\medskip We present an algorithm checking convexity of PL-realizations (in the sense outlined above) of a closed compact  manifold $M=({\cal M},r)$.

The main algorithm uses an auxiliary algorithm C-check. The input of this algorithm is
a pair $(T,{\cal T})$, where ${\cal T}$ is a one vertex tree with a cyclic orientation of edges and $T$ is its rectilinear realization  in 3-space. This pair can be thought of as a PL-realization of a plane fan (partition of the plane into cones with common origin)  in 3-space. The output is 1, if this realization is convex, and 0 otherwise. Obviously, this question is equivalent to verifying  convexity  of a plane polygon. For the plane of reference we choose a plane perpendicular to the sum of all unit vectors directed along the edges of the fan.  The latter question can be resolved in time, linear in the number of edges of the tree (e.g. see Devillers et al (1998), Mehlhorn et al (1999)).

\underline{Input and Preprocessing:} The poset of faces of dimensions $n-3,n-2$, and $n-1$  of ${\cal M}$ and the equations of the facets, OR the poset of faces of dimensions $n-3,n-2,n-1$, and $0$  of ${\cal M}$ and the positions of the vertices. We assume that we know the correspondence between the rank of a face in the poset and its dimension.  There are mutual links between the facets (or vertices) of ${\cal M}$ in the poset and the records containing their realization information. All $(n-3)$-faces of ${\cal M}$ are put into a stack $S_{n-3}$.
There are mutual links between elements of this stack and corresponding elements of the face lattice of ${\cal M}$.

\underline{Output:} YES, if $r({\cal M})$ is the boundary of a convex polyhedron, NO otherwise.

\begin{tabbing}

%\begin{itemize}

%\item
1. {\bf while} $S_{n-3}$ is not empty, pick an $(n-3)$-face $F$ from $S_{n-3}$;\\

%\item
2. compute the projection of $F$, and of all $(n-2)$-faces incident to $F$,\\
~~~onto an affine 3-plane complimentary to $F$; denote this 
projection  by $\PStar(F)$;\\

%\item
3. compute the cyclicly ordered one-vertex tree ${\cal T}(F)$, whose edges\\
~~~are the $(n-2)$-faces of $\PStar(F)$ and whose vertex is $F$; \\

%\item
4. Apply to $(\PStar(F),{\cal T}(F))$ the algorithm C-check\\

~~~~~{\bf if} C-check$(\PStar(F),{\cal T}(F))=1$ \= {\bf then} remove $F$ from the stack $S_{n-3}$ \\

\>    {\bf else} Output:=NO; terminate \\

%\item
~~~{\bf endwhile}\\

5. Output:=YES \\

%\end{itemize}

\end{tabbing}

\begin{remark} The algorithm processes the stars of all $(n-3)$-faces independently. On a parallelized computer the stars of all $(n-3)$-faces can be processed in parallel.
\end{remark}

\textbf{Proof of Correctness.} The algorithm checks the local convexity of $M$ at the stars of all $(n-3)$-cells.  $\cM$ is compact and closed — by Krein-Milman theorem (or Lemma \ref{strict}) $M$ has at least one strictly convex vertex.  By Theorem \ref{PL-case} local convexity at all vertices, together with the existence of at least one strictly convex vertex, is necessary and sufficient for $M$ to be the boundary of a convex body.

\underline{\bf Complexity estimates}
Denote by $f_k$ the number of $k$-faces of ${\cal M}$, and by $f_{k,l}$ -- the number of incidences between $k$-faces and $l$-faces in ${\cal }M$. Step 1 is repeated at most $f_{n-3}$ times. Steps 2-4 take  at most $\const f_{n-2,n-3} (\Star(F))$ arithmetic operations for each $F$, where $\const$ does not depend on $F$. Thus, steps 2-4, repeated for all $(n-3)$-faces of ${\cal M}$, require $O(f_{n-2,n-3})$ operations. Therefore, the total number of operations for this algorithm is $O(f_{n-2,n-3})$.

\begin{remark} The algorithm does not use all of the face lattice of ${\cal M}$. 
\end{remark}
\begin{remark} The algorithm requires computing polynomial predicates only. The highest degree of algebraic predicates that the algorithm uses is $d$, which is optimal (see Devillers et al, 1998). 
\end{remark}

From a practical point of view, it makes sense to say that a surface $M$ is almost convex, if it lies within a small Hausdorf distance from a convex surface $S$ that bounds an $n$-dimensional convex set $B$. In this case, the measure of lines, that pass through interior points of $B$ and intersect $S$ in more than 2 points, will be small, as compared to the measure of all lines passing through interior points of $B$. These statements can be given a rigorous meaning in the language of integral geometry, also called ``geometric probability'' (see Klain and Rota 1997).   

\begin{remark} If there is a 3-dimensional coordinate subspace $L$ of $\R^n$ such that all the subspaces spanned by $(n-3)$-faces are complementary to $L$, the polyhedron can be projected on $L$ and all computations can be done in 3-space. This reduces the degree of predicates from $d$ to 3. In such case the boolean complexity of the algorithm does not depend on the dimension at all. therefore, for sufficiently generic realizations the algorithm has degree 3 and complexity not depending on $n$.
\end{remark}
\begin{remark} This algorithm can also be applied without changes to compact PL-surfaces in $\SSS^n$ or $\HH^n$.
\end{remark}

\end{document}